\documentclass{article}
\PassOptionsToPackage{numbers}{natbib}
\usepackage{optml}
\usepackage{amsmath}

\usepackage[utf8]{inputenc} 
\usepackage[T1]{fontenc}    
\usepackage{hyperref}       
\usepackage{url}            
\usepackage{booktabs}       
\usepackage{amsfonts}       
\usepackage{nicefrac}       
\usepackage{microtype}      

\usepackage{amsmath,amsfonts,amssymb,amsthm,array}
\usepackage{amsmath}
\usepackage{algorithm}
\usepackage{caption}
\usepackage{subcaption}
\usepackage{algpseudocode}
\usepackage{bm}
\usepackage{color}
\usepackage{graphicx}

\newcommand{\Exp}{\mathbb{E}}

\newcommand{\R}{\mathbb{R}}

\usepackage{tikz}

\newcommand{\bA}{\mathbf{A}}
\newcommand{\bB}{\mathbf{B}}

\newcommand{\bS}{\mathbf{S}}

\newcommand{\eqdef}{:=}

\newcommand{\cD}{{\cal D}}

\newcommand{\cO}{{\cal O}}

\newcommand{\mA}{{\bf A}}

\newcommand{\mH}{{\bf H}}

\newcommand{\mS}{{\bf S}}

\theoremstyle{plain}
\newtheorem{thm}{Theorem}[]

\newtheorem{rem}{Remark}[]

\theoremstyle{remark}

\usepackage[colorinlistoftodos,bordercolor=orange,backgroundcolor=orange!20,linecolor=orange,textsize=scriptsize]{todonotes}

\usepackage{hyperref}

\title{Linearly Convergent Stochastic Heavy Ball Method for Minimizing Generalization Error}
\author{\name Nicolas Loizou \email{n.loizou@sms.ed.ac.uk}\\
  \addr{ University of Edinburgh, United Kingdom}\\
  \name Peter Richt\'{a}rik  \email{peter.richtarik@ed.ac.uk}\\
  \addr{ KAUST, Kingdom of Saudi Arabia \\
  University of Edinburgh, United Kingdom}
}

\begin{document}
\maketitle

\begin{abstract}
 In this work we establish the first linear convergence result for the stochastic heavy ball method. The method performs SGD steps with a fixed  stepsize, amended by a heavy ball momentum term. In the analysis, we focus on minimizing the expected  loss and not on  finite-sum minimization, which is typically a much harder problem. While in the analysis we constrain ourselves to quadratic loss, the overall objective is not necessarily strongly convex.
\end{abstract}

\section{Introduction}
In this paper we study the stochastic optimization problem:
\begin{equation}
\label{eq:stoch_reform}
\min_{x\in \R^d} f(x) \eqdef \Exp_{\mS\sim \cD}[f_\mS(x)]
\end{equation}
where $\mA\in \R^{m\times d}$ is a data matrix, $b\in \R^m$ is a vector of labels, $\mS$ is a matrix with $m$ rows (and arbitrary number of columns, e.g., 1), $\cD$ is a distribution over such matrices and $f_{\mS}(x) \eqdef \tfrac{1}{2}\|\mA x - b\|_{\mH}^2$  is a least-squares function with respect to a random pseudo-norm defined by a specific symmetric positive semidefinite  matrix $\mH$ which depends on $\mA$ and the random matrix $\mS$.  In particular, $\|y\|_{\mH}^2 \eqdef y^\top \mH y$ and $\mH \eqdef  \mS (\mS^\top \mA \mA^\top \mS)^\dagger \mS^\top$, where $\dagger$ denotes the Moore-Penrose pseudoinverse. Note that the function $f$ is finite if and only if $\Exp_{\mS\sim \cD}[\mH]$ exists and is finite. Hence, we assume this throughout this paper. 

Problem \eqref{eq:stoch_reform} was first proposed in \cite{ASDA}, where the authors focus on stochastic reformulations of a consistent linear system $\mA x = b$. The authors further give necessary and sufficient conditions on $\cD$ for the set of solutions of \eqref{eq:stoch_reform} to be equal to the set of solutions of the linear system $\mA x = b$; a property for which the term \emph{exactness} was coined in \cite{ASDA}. Exactness conditions are very weak, allowing $\cD$ to be virtually any distribution of random matrices. For instance, a sufficient condition for exactness is for the matrix $\Exp[\mH]$ to be positive definite. This is indeed a weak condition since it is easy to see that  this matrix is symmetric and positive semidefinite without the need to invoke any assumptions; simply by design. We refer the reader to \cite{ASDA} for more insights into the reformulation \eqref{eq:stoch_reform}, its properties and other equivalent reformulations (e.g., stochastic fixed point problem, probabilistic intersection problem, and stochastic linear system). 

In \cite{ASDA}, the authors  consider solving \eqref{eq:stoch_reform} via stochastic gradient descent (SGD)
\begin{equation}\label{eq:SGD}x_{k+1} = x_k - \omega \nabla f_{\mS_k}(x_k),\end{equation}
where $\omega>0$ is a fixed stepsize and  $\mS_k$ is sampled afresh in each iteration from $\cD$. It is shown that, SGD converges to an $x_*$ which satisfies \begin{equation}\label{eq:primal}x_*=\text{argmin}_{x\in \R^d}  \tfrac{1}{2}\|x-x_0\|^2 \quad \text{subject to} \quad \mA x = b ,\end{equation} where $x_0$ is the starting point. It was observed that, surprisingly, SGD is in this setting equivalent to the stochastic (pseudo)-Newton method, and the stochastic proximal point method, and that it converges at a linear rate despite the following obstacles: $f$ is not necessarily strongly convex,  \eqref{eq:stoch_reform} is not a finite-sum problem, and a fixed stepsize $\omega$ is used.

\subsection{Contributions}
In this paper we take an alternative route, and develop a {\em stochastic} variant of the {\em heavy ball method} for solving the stochastic optimization problem  \eqref{eq:stoch_reform}. Applied to  \eqref{eq:stoch_reform}, the classical heavy ball method of Polyak \cite{polyak1964some, polyak1987introduction}, with constant stepsize $\omega>0$ and constant momentum parameter $\beta\geq 0$, takes the form
\begin{equation}
\label{HB}
x_{k+1} = x_k - \omega \nabla f(x_k) + \beta(x_k - x_{k-1})
\end{equation}
This method  introduces the momentum term  $\beta(x_k - x_{k-1})$ into the gradient descent method to achieve acceleration.

Our stochastic variant of the heavy ball method, which we henceforth simply refer to by the name {\em stochastic heavy ball method (SHB)}, replaces the (costly) computation of the gradient by an unbiased estimator of the gradient (``stochastic gradient'') which is hopefully much cheaper to compute:
\begin{equation}\label{eq:SHB-intro} \boxed{x_{k+1} = x_k - \omega \nabla f_{\mS_k}(x_k) + \beta(x_k-x_{k-1})}\end{equation}
We establish global linear convergence (in expectation) of the iterates and function values: $\Exp[\|x_k-x_*\|^2] \to 0$ (L2 convergence) and $\Exp[f(x_k)] \to 0$. Without the exactness assumption we prove that  $\Exp[f(\hat{x}_k)]=\cO(1/k)$, where  $\hat{x}_k= \tfrac{1}{k}\sum_{t=0}^{k-1} x_t$ is the Cesaro average. Finally, we study the convergence of the expected iterates (L1 convergence), $\|\Exp[x_k-x_*]\|^2\to 0$,  and establish global accelerated linear rate. That is, this quantity falls below $\epsilon$ after $\cO\left((\lambda_{\max}/\lambda_{\min}^+)^{1/2} \log(1/\epsilon)\right)$ iterations, where $\lambda_{\max}$ (resp.\ $\lambda_{\min}^+$) are the largest (resp.\ smallest nonzero) eigenvalues of: $\nabla^2 f(x) =  \mA^\top \Exp_{\mS\sim \cD}[\mH] \mA.$ It turns out that all eigenvalues of $\nabla^2 f(x)$ belong to the interval $[0,1]$.

\subsection{Related Work}
Stochastic variants of heavy ball method have been employed widely in practice, especially in the area of deep learning \cite{sutskever2013importance, szegedy2015going, krizhevsky2012imagenet}. Despite the popularity of the method both in convex and non-convex optimization its convergence properties are not very well understood.  Recent papers that provide complexity analysis of SHB (in different setting than ours) include \cite{yang2016unified} and \cite{gadat2016stochastic}. In \cite{yang2016unified} the authors analyzed  SHB for general Lipshitz continuous convex objective functions (with bounded variance) and proved the {\em sublinear} rate $O(1/\sqrt{k})$. In \cite{gadat2016stochastic}, a complexity analysis is provided for the case of quadratic strongly convex smooth coercive functions. A {\em sublinear} convergence  rate $O(1/k^\beta)$, where $\beta \in (0,1)$, was proved. In contrast to our results, where we assume fixed stepsize $\omega$, both papers analyze SHB with diminishing stepsizes. For our problem, variance reduction methods like SVRG \cite{johnson2013accelerating}, S2GD \cite{S2GD}, mS2GD \cite{mS2GD}, SAG \cite{schmidt2017minimizing} and SAGA \cite{defazio2014saga} are not necessary. To the best of our knowledge, {\em our work provides the first linear convergence rate for SHB in any setting.}

\section{Convergence Results}

\label{sec:SHB}
In this section we state our convergence results for SHB. 

\subsection{$L2$ convergence: linear rate}
We study L2 convergence of SHB; that is, we study the convergence of the quantity $\Exp[\|x_k-x_*\|^2]$ to zero. We show that for a range  of stepsize parameters $\omega > 0$ and  momentum parameters $\beta \geq 0$, SHB enjoys \emph{global non-asymptotic linear convergence rate}. As a corollary of  L2 convergence, we obtain convergence of the expected function values.
\begin{thm} 
\label{L2}
Choose $x_0= x_1\in \R^d$.  Assume exactness. Let $\{x_k\}_{k=0}^\infty$ be the sequence of random iterates produced by SHB.  Assume $0< \omega < 2$ and $\beta \geq 0$ and that the expressions
\[a_1 \eqdef 1+3\beta+2\beta^2 - (\omega(2-\omega) +\omega\beta)\lambda_{\min}^+, \qquad \text{and}\qquad
a_2 \eqdef \beta +2\beta^2 + \omega \beta \lambda_{\max}\]
satisfy $a_1+a_2<1$. Let $x_*$ be the projection of $x_0$ onto $\{x\;:\; \bA x = b\}$. Then 
\begin{equation}\label{eq:nfiug582}\Exp[\|x_{k}-x_*\|^2] \leq q^k (1+\delta)  \|x_{0}-x_*\|^2\end{equation}
and 
$$\Exp[f(x_k)] \leq q^k  \tfrac{\lambda_{\max}}{2} (1+\delta) \|x_{0}-x_*\|^2,$$
where  $q=\frac{a_1+\sqrt{a_1^2+4a_2}}{2}$ and $\delta=q-a_1$. Moreover, $a_1+a_2 \leq q <1$.
\end{thm}

\begin{rem}
In the above theorem we obtain global linear rate. To the best of our knowledge, this is the first time that linear rate is established for a stochastic variant of the heavy ball method in any setting. All existing results are sublinear.
\end{rem}

\begin{rem}
\label{rangesSHB}
If we choose $\omega \in (0,2)$, then the condition $a_1+a_2<1$ is satisfied for all   $$0\leq \beta< \tfrac{1}{8} \left( -4+\omega \lambda_{\min}^+-\omega \lambda_{\max} +\sqrt{(4-\omega \lambda_{\min}^++\omega \lambda_{\max})^2+16\omega (2-\omega) \lambda_{\min}^+ }\right).$$
\end{rem}

\begin{rem}
If $\beta=0$, SHB reduces to the ``basic method'' in \cite{ASDA} (SGD with constant stepsize). In this special case, $q = 1-\omega(2-\omega)\lambda_{\min}^+$, which is the rate established in \cite{ASDA}. Hence, our result is more general.
\end{rem}

\begin{rem} Let $q(\beta)$ be the rate as a function of $\beta$. Note that since $\beta\geq 0$, we have
\begin{equation}\label{eq:qbeta} q(\beta) \geq a_1 + a_2 = 1 + 4\beta + 4\beta^2 + \omega\beta(\lambda_{\max}-\lambda_{\min}^+) - \omega(2-\omega)\lambda_{\min}^+ \geq 1-\omega(2-\omega)\lambda_{\min}^+ = q(0).\end{equation}
Clearly, the lower bound on $q$ is an increasing function of $\beta$. 
Also, for any $\beta$ the rate is always inferior to that of SGD ($\beta=0$). It is an open problem whether one can prove a strictly better rate for SHB than for SGD.
\end{rem}

\subsection{Cesaro average: sublinear rate without exactness assumption}
In this section we present convergence results for function values computed at the Cesaro average of all past iterates. Again, our results are global in nature. To the best of our knowledge, an analysis of the Cesaro average for the SHB with $O(1/k)$ rate was not established before for any class of functions. Moreover, the result holds without the exactness assumption.
\begin{thm}
\label{cesaro}
Choose $x_0=x_1$ and let $\{x_k\}_{k=0}^\infty$ be the random iterates produced by SHB, where the momentum parameter $0\leq \beta <1$ and relaxation parameter (stepsize) $\omega > 0$ satisfy $\omega + 2\beta <2$. Let $x_*$ be any vector satisfying $f(x_*)=0$. If we let $\hat{x}_k=\frac{1}{k}\sum_{t=1}^{k}x_t$, then
$$\Exp[f(\hat{x}_k)] \leq \frac{(1-\beta)^2\|x_0-x_*\|^2 + 2\omega \beta f(x_0)}{2\omega(2-2\beta-\omega) k}.$$
\end{thm}
\begin{rem}
In the special case of $\beta=0$ we have  
$\Exp[f(\hat{x}_k)] \leq \tfrac{\|x_0-x_*\|^2}{2\omega (2 -\omega) k}$,
which is the convergence rate for Cesaro average of the ``basic method'' analyzed in \cite{ASDA}. 
\end{rem}
\subsection{$L1$ convergence: accelerated linear rate}
In this section we show that by a proper combination of the stepsize parameter $\omega$ and the momentum parameter $\beta$ the proposed algorithm enjoys {\em accelerated linear convergence} rate with respect to the expected iterates.
\begin{thm}
\label{theoremheavyball}
Assume exactness. Let $\{x_k\}_{k=0}^{\infty}$ be the sequence of random iterates produced SHB, started with $x_0, x_1 \in \R^d$ satisfying the relation $x_0-x_1 \in {\rm Range}(\bA^ \top)$, with stepsize parameter  $0<\omega \leq1/\lambda_{\max}$ and momentum parameter  $\left(1-(\omega \lambda_{\min}^+)^{1/2}\right)^2 < \beta <1$. Then there exists constant $C >0$ such that for all $k\geq0$ we have 
$\|\Exp[x_{k} -x_*]\|^2  \leq \beta^k C.$
\begin{itemize}
\item[(i)] If we choose $ \omega= 1$ and $\beta= \left(1- \sqrt{0.99 \lambda_{\min}^+}\right) ^2$, then
 $$\|\Exp[x_{k} -x_*]\|_{\bB}^2  \leq \left(1- \sqrt{0.99 \lambda_{\min}^+}\right) ^{2k} C$$
and the iteration complexity becomes
$ \cO\left(\sqrt{1/ \lambda_{\min}^+}\log(1/\epsilon)\right)$.
\item[(ii)] If we choose $ \omega= 1/\lambda_{\max}$ and $\beta= \left(1- \sqrt{ 0.99 \lambda_{\min}^+/\lambda_{\max}}\right)^2$, then
$$\|\Exp[x_{k} -x_*]\|_{\bB}^2  \leq \left(1- \sqrt{ 0.99 \lambda_{\min}^+/\lambda_{\max}}\right)^{2k} C$$
and the iteration complexity becomes $\cO\left(\sqrt{\lambda_{\max}/ \lambda_{\min}^+} \log(1/\epsilon)\right)$
\end{itemize}
\end{thm}
Note that the convergence factor is precisely equal to the value of the momentum parameter.

\begin{rem} Let $x$ be any random vector in $\R^d$ with finite mean $\Exp[x]$, and $x_*\in \R^d$ be any reference vector (for instance, any solution of $\mA x = b$). Then we have the identity (see, for instance \cite{gower2015randomized})
$$\Exp[\|x-x_*\|^2] = \|\Exp[x-x_*]\|^2 + \Exp[\|x-\Exp[x]\|^2].$$
This means that the quantity $\Exp[\|x-x_*\|^2]$ appearing in our L2 convergence result (Theorem~\ref{L2}) is larger than 
$\|\Exp[x-x_*]\|^2$ appearing in the L1 convergence result (Theorem~\ref{theoremheavyball}), and hence harder to push to zero. As a corollary, L2 convergence implies L1 convergence. However, note that in Theorem~\ref{theoremheavyball} we have established an {\em accelerated} rate. 
\end{rem}

\section{Experiments}
\begin{figure}[H]
\centering
\begin{subfigure}{.23\textwidth}
  \centering
  \includegraphics[width=1\linewidth]{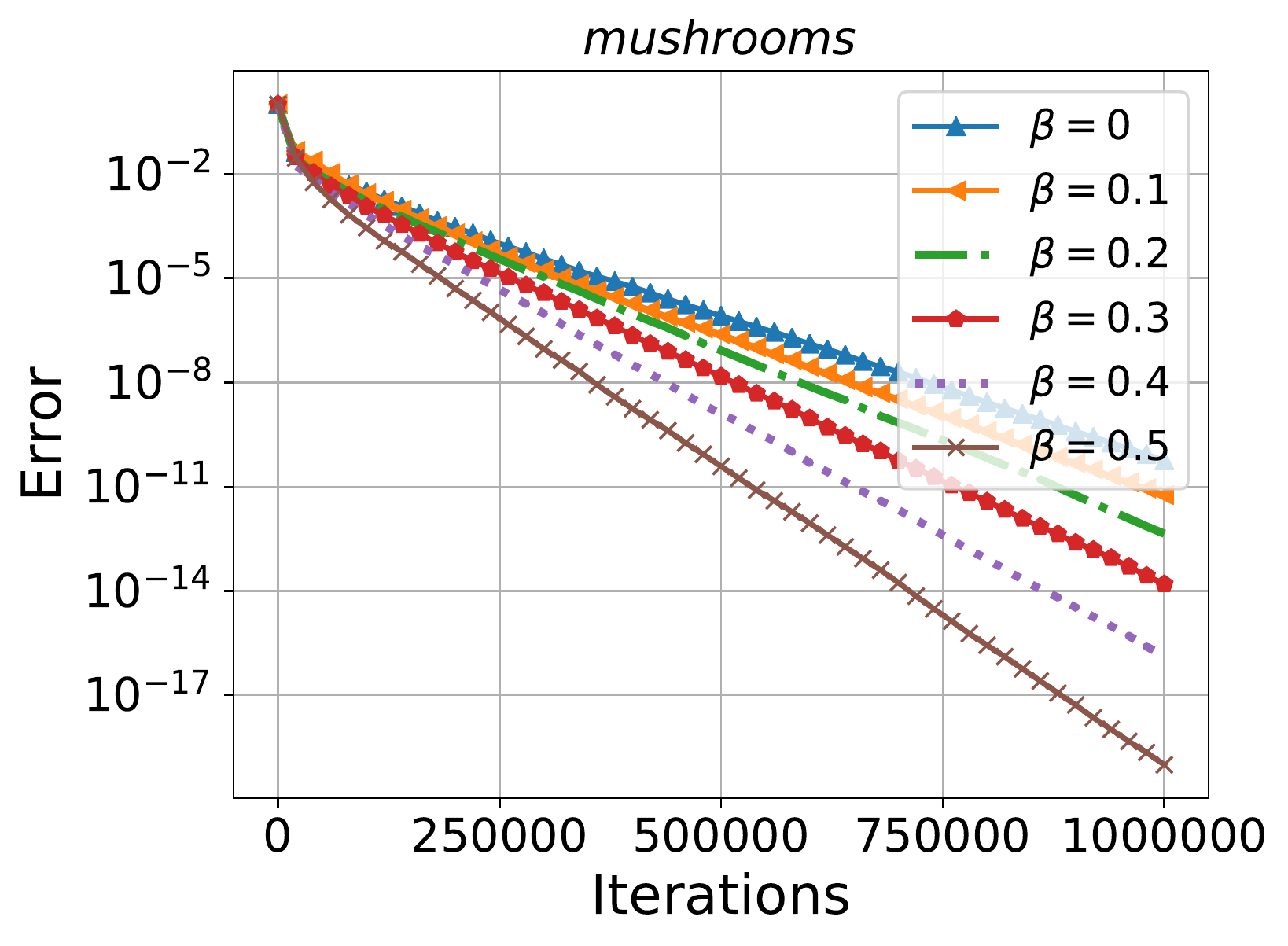}
\end{subfigure}%
\begin{subfigure}{.23\textwidth}
  \centering
  \includegraphics[width=1\linewidth]{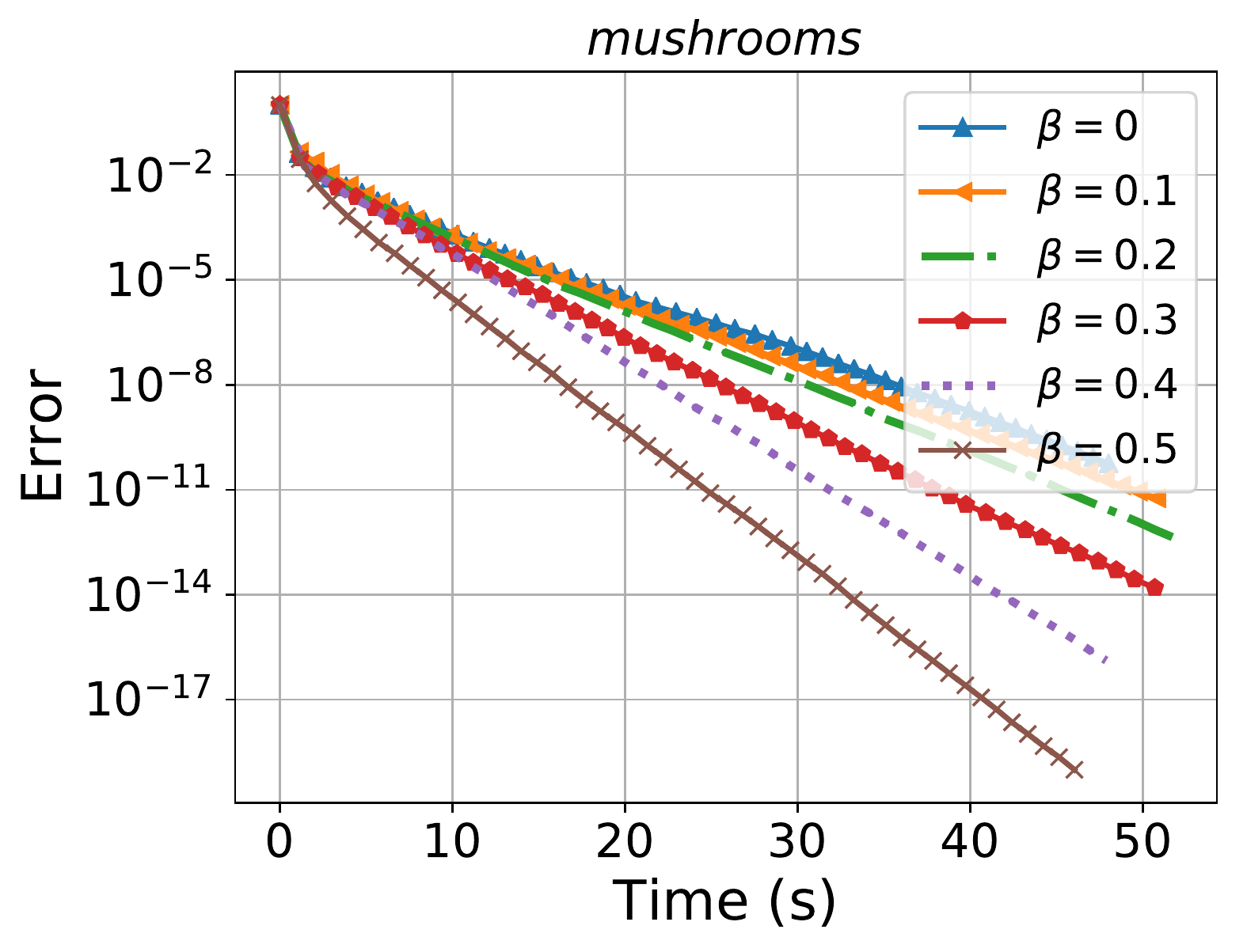}
\end{subfigure}
\begin{subfigure}{.23\textwidth}
  \centering
  \includegraphics[width=1\linewidth]{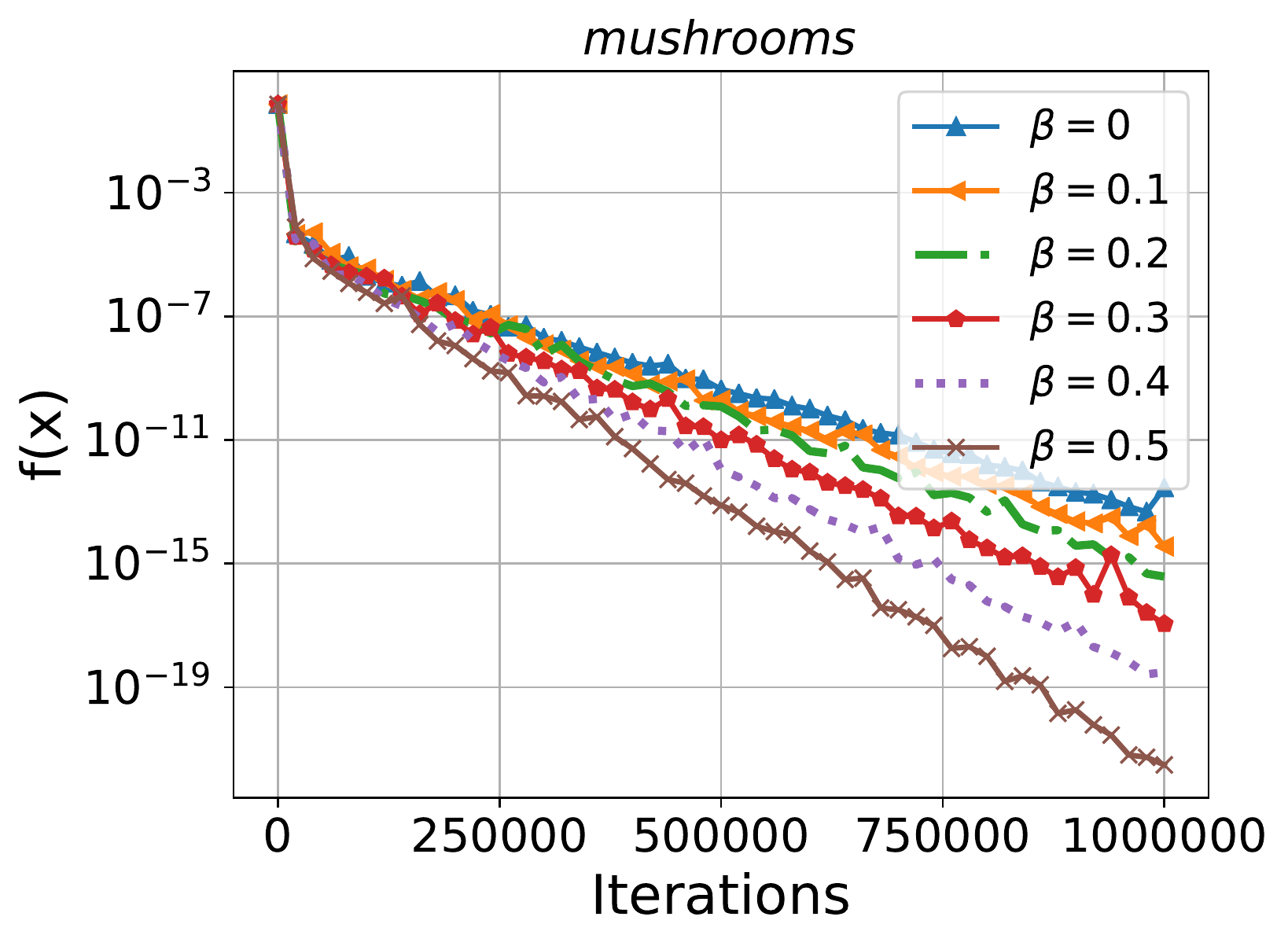}
\end{subfigure}
\begin{subfigure}{.23\textwidth}
  \centering
  \includegraphics[width=1\linewidth]{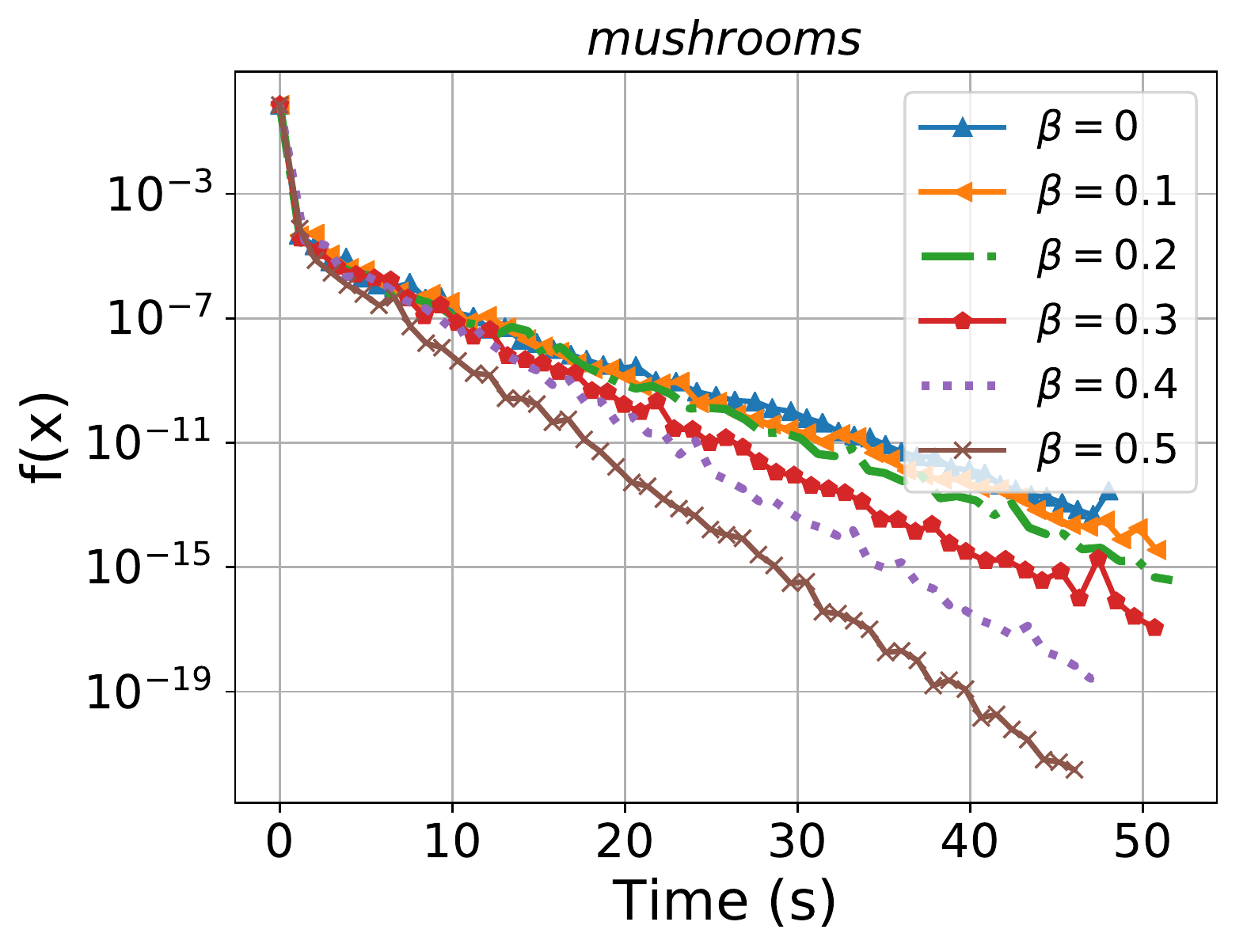}
\end{subfigure}\\
\begin{subfigure}{.23\textwidth}
  \centering
  \includegraphics[width=1\linewidth]{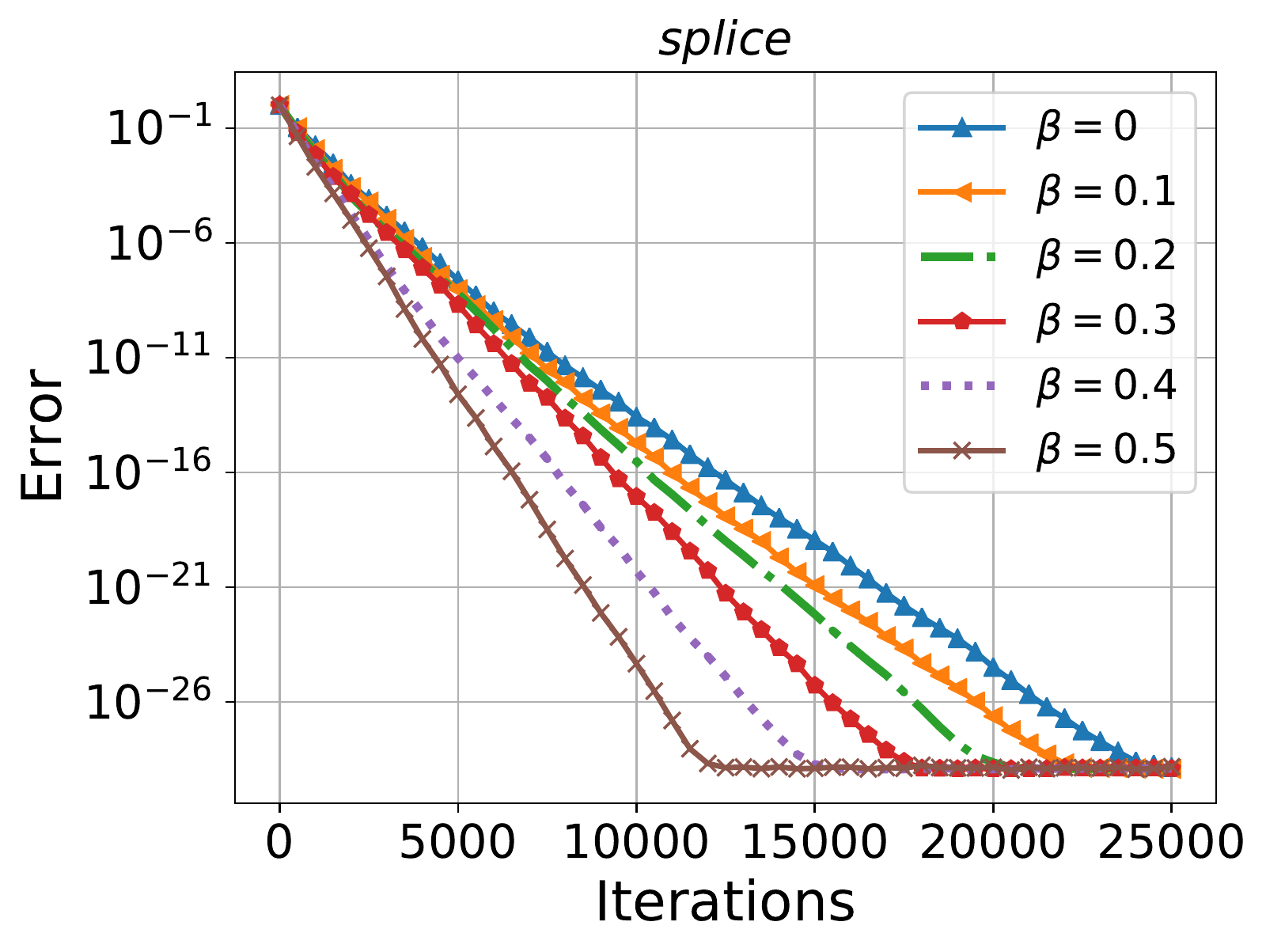}
\end{subfigure}%
\begin{subfigure}{.23\textwidth}
  \centering
  \includegraphics[width=1\linewidth]{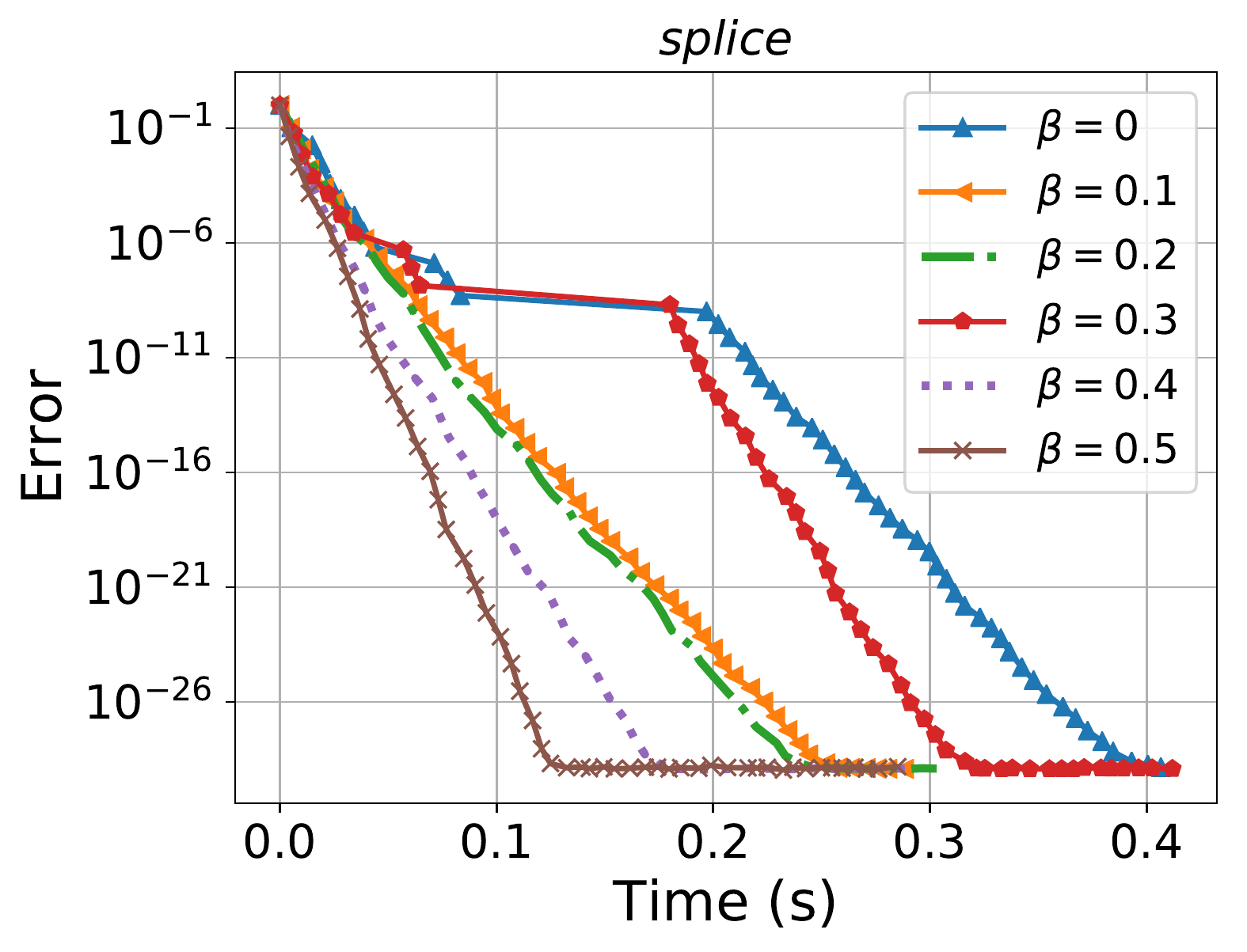}
\end{subfigure}
\begin{subfigure}{.23\textwidth}
  \centering
  \includegraphics[width=1\linewidth]{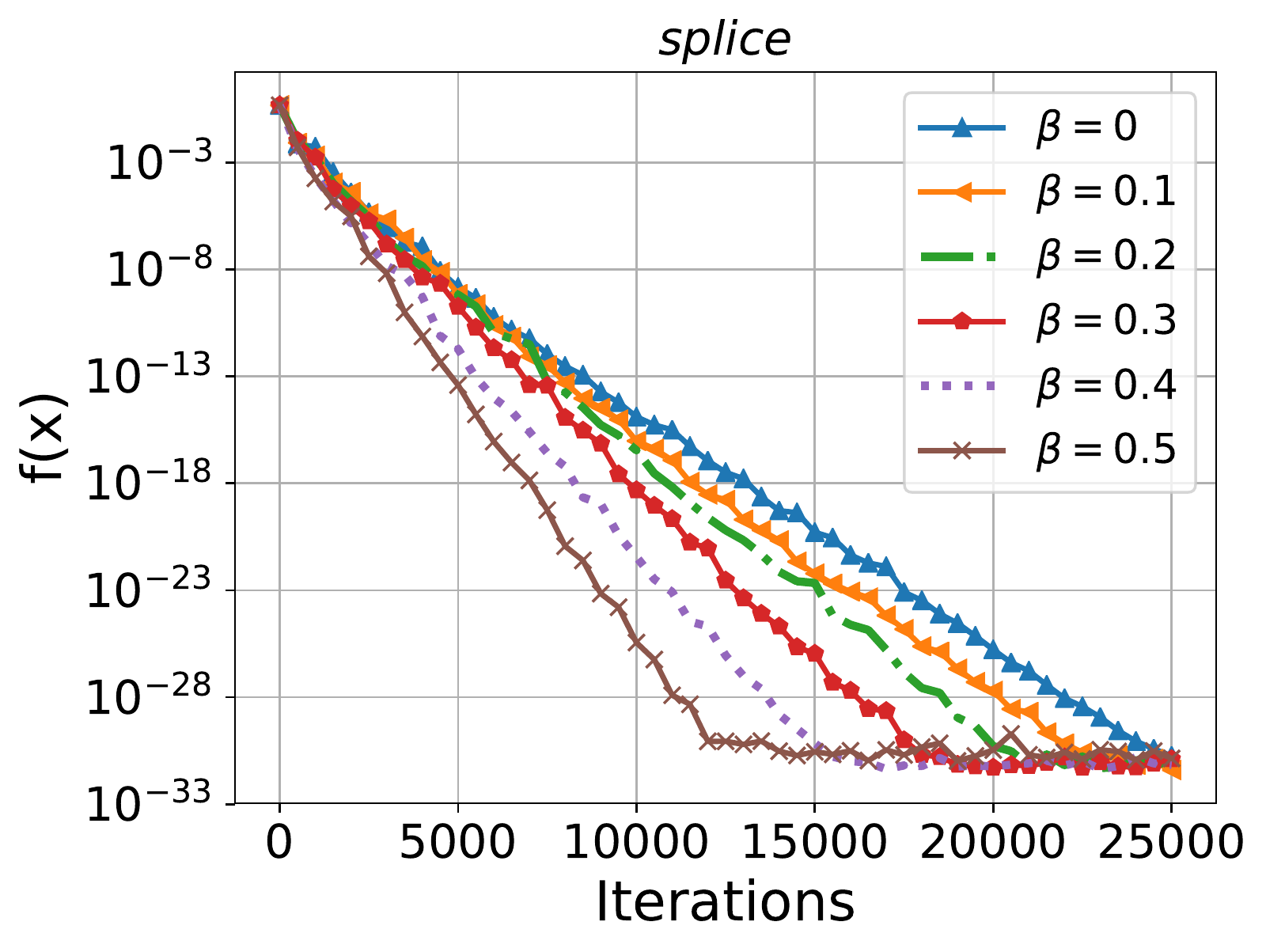}
\end{subfigure}
\begin{subfigure}{.23\textwidth}
  \centering
  \includegraphics[width=1\linewidth]{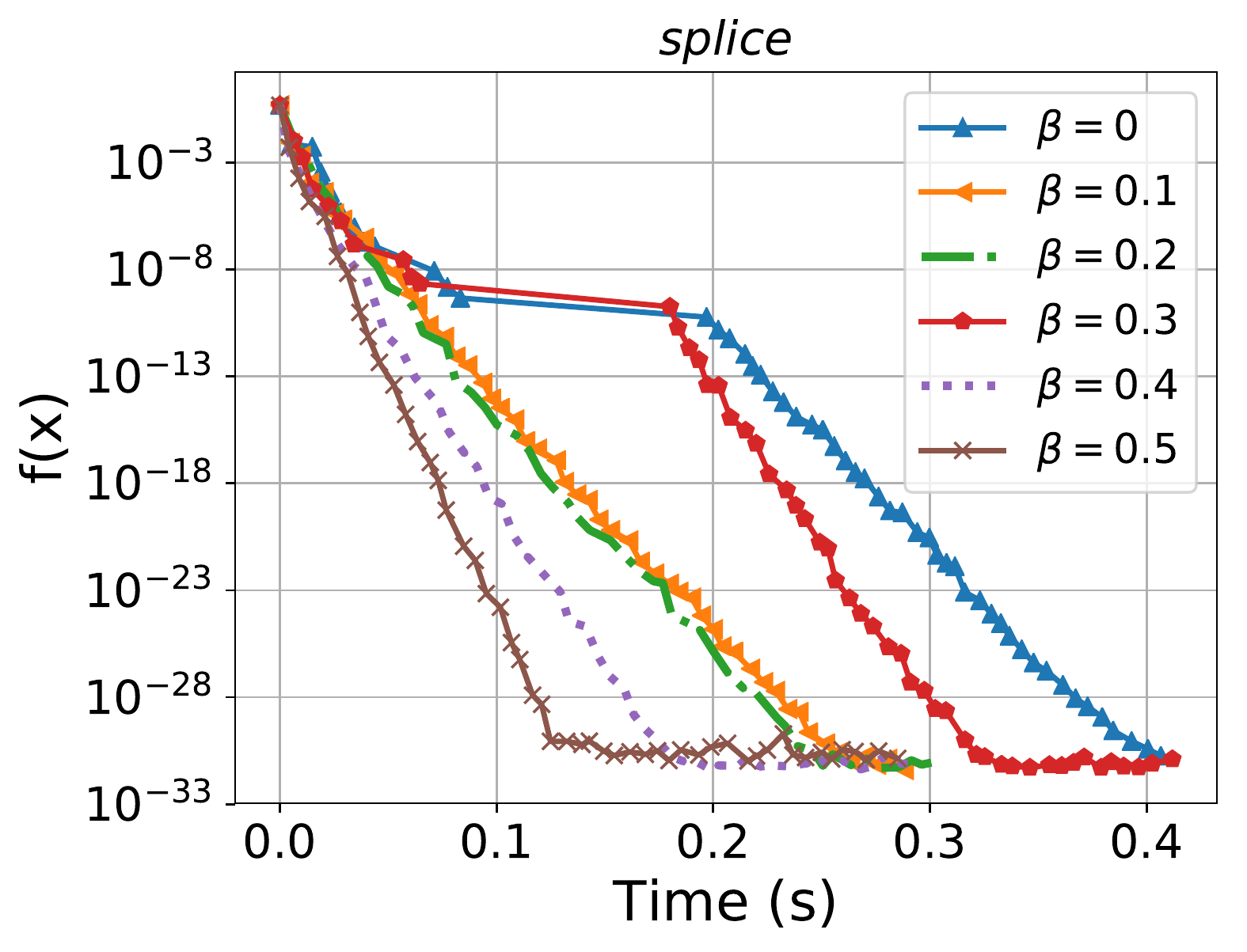}
\end{subfigure}\\
\caption{The performance of RK with momentum for several momentum parameters $\beta$ on real data from LIBSVM \cite{chang2011libsvm}. mushrooms: $(n,d)=(8124,112)$, splice: $(n,d)=(1000,60)$. The graphs in the first (second) column plot iterations (time) against residual error while those in the third (forth) column plot iterations (time) against function values. The ``Error" on the vertical axis represents the relative error $\frac{\|x_k-x_*\|^2}{\|x_*\|^2}$ and the function values $f(x_k)$ refer to function \eqref{eq:stoch_reform}.}
\label{RealDataplots}
\end{figure}

In this section we present a preliminary experiment to evaluate the performance of the SHB for solving  the stochastic optimization problem~\eqref{eq:stoch_reform}. Matrices $\bA$ are picked from the LIBSVM library \cite{chang2011libsvm}. To ensure consistency of the linear system, we take the optimal solution $x_{*}\in \R^d$ to be i.i.d $\mathcal{N}(0,1)$ and the right hand sight is set to $b=\bA x_{*}$. 
In each iteration, the random matrix is chosen as $\bS=e_i \in \R^n$ with probability $p_i=\|\bA_{i:}\|^2 / \|\bA\|_F^2$. Here $e_i$ is the unit coordinate vector in $\R^{n}$. In this setup the update rule \eqref{eq:SHB-intro} of the SHB simplifies to \[x_{k+1}=x_k -\omega \frac{\bA_{i:} x_k -b_i}{\|\bA_{i:}\|_2^2} \bA_{i:}^ \top + \beta(x_k - x_{k-1}) .\] This is a randomized Kaczmarz method (RK) with momentum. Note that for $\beta=0$ and $ \omega=1$ this reduces to the celebrated \textit{Randomized Kaczmarz method} (RK) of Strohmer and Vershynin \cite{RK}. In Figure~\ref{RealDataplots}, RK with momentum is tested for several values of the momentum parameters $\beta$ and fixed stepsize $\omega=1$.  For the evaluation we use both the relative error measure $\|x_k-x_*\|^2/ \|x_0-x_*\|^2$ and the function suboptimality $f(x_k)-f(x_*)$. The starting point is chosen as $x_0=0$. For the horizontal axis we use either the number of iterations or the wall-clock time measured using the tic-toc Julia function. It is clear that in this setting the addition of momentum parameter is beneficial and leads to faster convergence. 
\newpage
\small
\bibliographystyle{plain}
\bibliography{biblio}

\end{document}